# Mathematical Modeling and System Identification Aspects of Heat Flow Process inside a Closed Box

Sulaymon Eshkabilov, Elyor Madirimov
Dynamics & Control Lab
Tashkent Automotive Road Institute
Amir Temur Str. 20, Tashkent -100060, Uzbekistan
Email: sulaymon@d-c-lab.com

***Abstract.*** *This paper presents review points of mathematical modeling and system identification issues of heat flow process inside a closed box. For system parameter identification, the least squares method is employed and NI® Temperature Box is used as a test rig for different reference temperature values. Curve Fitting™ toolbox and M-file codes written in MATLAB are used based on the Levenberg-Marquardt algorithm in solving least squares equations of the system.*

## 1. Introduction

Mathematical model or formulation of given system variables is very important not only to understand the system behaviors but also to control it efficiently and predict its outputs or its responses on reference input signals. Thus, we strive to identify the relationship between system variables and parameters when we start designing new systems or making changes on existing ones. To achieve this objective is not straightforward and there is no a unique algorithm how to reach it for all systems. There are different approaches and methods for system modeling and identification depending on the system type or in other words how much knowledge we have about the system and how much noiseless reference input-output data we have at our disposal. Also, if the given system is white box, grey box or black box model type. How we are intending to use our desired model will define its complexity [1].

In this paper, we and derive general mathematical formulations of a heat flow process inside a closed box using experimental test measurements of input-out data in real time and define system parameters by performing curve fitting with the least squares method and using the Levenberg-Marquardt algorithm in solving the equations of the least squares method. We have considered NI® Temperature box that is a laboratory test rig to demonstrate heat flow inside a closed box environment. Our chosen model is a grey-box model type that has some properties known and some unknown. In order to remove undesired noises in our measured data sets we have used ready to use low-pass finite impulse Savitky-Golay filter of MATLAB package.

In fact, heat exchangers are complex systems that involve nonlinearities [2-4]. The heat flow process itself is non-linear that will be modeled similar to many other real time experimental data based modelling and require several steps to take. First step is to select modeling technique, second one is a model structure and the last one is model estimation [5]. The work [6] is dedicated to modeling and simulating a heating system composed of halogen lamp, energy source (current/voltage) and measuring unit (thermocouple) by employing the system identification toolbox of MATLAB. Another work [7] is devoted to modeling and system identification of home refrigerators using maximum likelihood estimation of parameters to design smart appliances by employing grey-box modeling approaches. The authors of the work [8] address the issue of modeling several types of heat exchangers by using ARX (autoregressive with exogenous input) model structure method in MATLAB environment.

## 2. Modeling a heat flow process

In our study we use thermal energy conceptions therefore we need to know basic law of thermal energy and it is known as the first law of thermodynamics. The first law of thermodynamics is also known as the conservation of energy principle, states that energy neither can be created nor destroyed; it can only change forms. Therefore, every bit of energy must be accounted for during a process. The conservation of energy principle (or energy balance) for any system undergoing any process can be expressed by the following expression:

| *Total energy entering the system* | - | *Total energy leaving the system* | = | *Change in the total energy of the system* |
|---|---|---|---|---|

Noting that energy can be transferred to or from a system by heat, work and mass flow. The energy balance for our system undergoing heating process can be expressed by

$$Q_{in} - Q_{out} = \Delta Q_{system} \tag{1}$$

or

$$Q_{gen} - Q_{loss} = Q_s \tag{2}$$

Where $Q_{gen}$ - energy generated by lamp, [J]; $Q_{loss}$ - energy loss to surrounding, [J]; $Q_s$ –stored energy in the box, [J].

***Heat source.*** Heating source in our system is a halogen lamp that gives fixed the temperature into a box. For calculating heat generation power in [W]:

$$\dot{Q}_{gen} = k \cdot V(t) \tag{3}$$

Where $k$ is a halogen lamp constant; $V(t)$ is a halogen lamp voltage.

***Heat loss.*** Let's consider steady-state heat transfer in a plane wall of thickness $\Delta x$. The wall temperature $T_1$ on one side (ambient temperature) and $T_2$ on the other (inside temperature of the temperature box). If the thermal conductivity is considered to be a constant k, then we substitute Fourier's law and result in [W]:

$$\dot{Q}_{loss} = A \cdot U \cdot (T - T_a) \tag{4}$$

Where $A$ is a surface area though which heat is transferred or lost from indoor to outdoor, [m$^2$]; $U$ is overall heat transfer coefficient of the material, in our case, material is Plexiglas, [W/m$^2$]; $T_a$ is a constant ambient temperature, [K]; $T$ is a indoor temperature of the box, which is changing by time, [K].

Stored heat in the system in [W] is expressed by the following formulation

$$\dot{Q}_s = \frac{dE}{dt} = \rho \cdot C_p \cdot \frac{dT}{dt} \tag{5}$$

Where $\rho$ is density of air, [kg/m$^3$]; $C_p$ is thermal capacity of air, $[J/kg \cdot K]$.

Finally, by substituting the expression (2) with the expressions (3), (4) and (5), we obtain our system's mathematical model

$$k \cdot V(t) = A \cdot U \cdot (T(t) - T_a) + \rho \cdot C_p \cdot \frac{dT}{dt} \tag{6}$$

From the first-order differential equation in (6) with "zero" initial conditions, we can switch it from time domain to "s" domain by applying Laplace transforms. We can obtain the transfer function of the system

$$k \cdot V(s) = A \cdot U \cdot T(s) + \rho \cdot C_p \cdot s \cdot T(s) \tag{7}$$

or

$$\frac{T(s)}{V(s)} = \frac{k/(A \cdot U)}{\left(\frac{\rho \cdot C_p}{A \cdot U}\right) \cdot s + 1} \tag{8}$$

However, for our simulations we take non-zero initial conditions, for instance, the ambient temperature $T_a = 25$°C, and thus, the transfer function of the process $T(s)$ will be equal to

$$T(s) = \frac{\frac{k}{A \cdot U} \cdot V(s)}{\left(\frac{\rho \cdot Cp}{A \cdot U}\right) \cdot s + 1} + \frac{T_a}{\left(\frac{\rho \cdot Cp}{A \cdot U}\right) \cdot s + 1} \tag{9}$$

If we substitute constants with $K=\frac{k}{A\cdot U}$ and $\tau=\frac{\rho\cdot C_p}{A\cdot U}$, we get simplified transfer function

$$T(s)=\frac{K\,V(s)}{\tau\cdot s+1}+\frac{T_a}{\tau\cdot s+1}=\frac{T_a+K\,V(s)}{\tau\cdot s+1} \tag{10}$$

If there is a time delay $t_d$ (in seconds) in the system, our transfer function will have the following form

$$T(s)=\frac{K\,V(s)\,e^{-t_d\cdot s}}{\tau\cdot s+1}+\frac{T_a}{\tau\cdot s+1} \tag{11}$$

The system models in (10) and (11) expressed in continuous time domain can be also re-formulated in discrete time domain by using three different approximation methods, such as, Tustin (bilinear transformation), forward and backward rectangular approximation methods.

The transfer functions of the system (the process) given in (10) and (11) are system models in continuous time domain with and without time delays in the system. From the transfer function of T(s) expression in (10) for Heaviside step input $V(s)=1/s$, we can derive the system response

$$T(t)=\frac{T_a e^{-\frac{t}{\tau}}}{\tau}+K\left(1-e^{-\frac{t}{\tau}}\right)=\frac{T_a e^{-\frac{t}{\tau}}}{\tau}+K-Ke^{-\frac{t}{\tau}} \tag{12}$$

Now we work with the test bench and collect measurement data from it in order to identify the system parameters, viz. $T_a, K, \tau$ using the models in time by employing Curve Fitting toolbox™ and writing a code in MATLAB® based on the least squares method with the Levenberg-Marquardt algorithm.

## 3. Experimental test bench description

The test bench is NI® temperature box composed of plexi-glass box, halogen lamp, J type thermo-couple and NI® PXI system with NI® LabView package and display monitor as shown in Fig. 1.

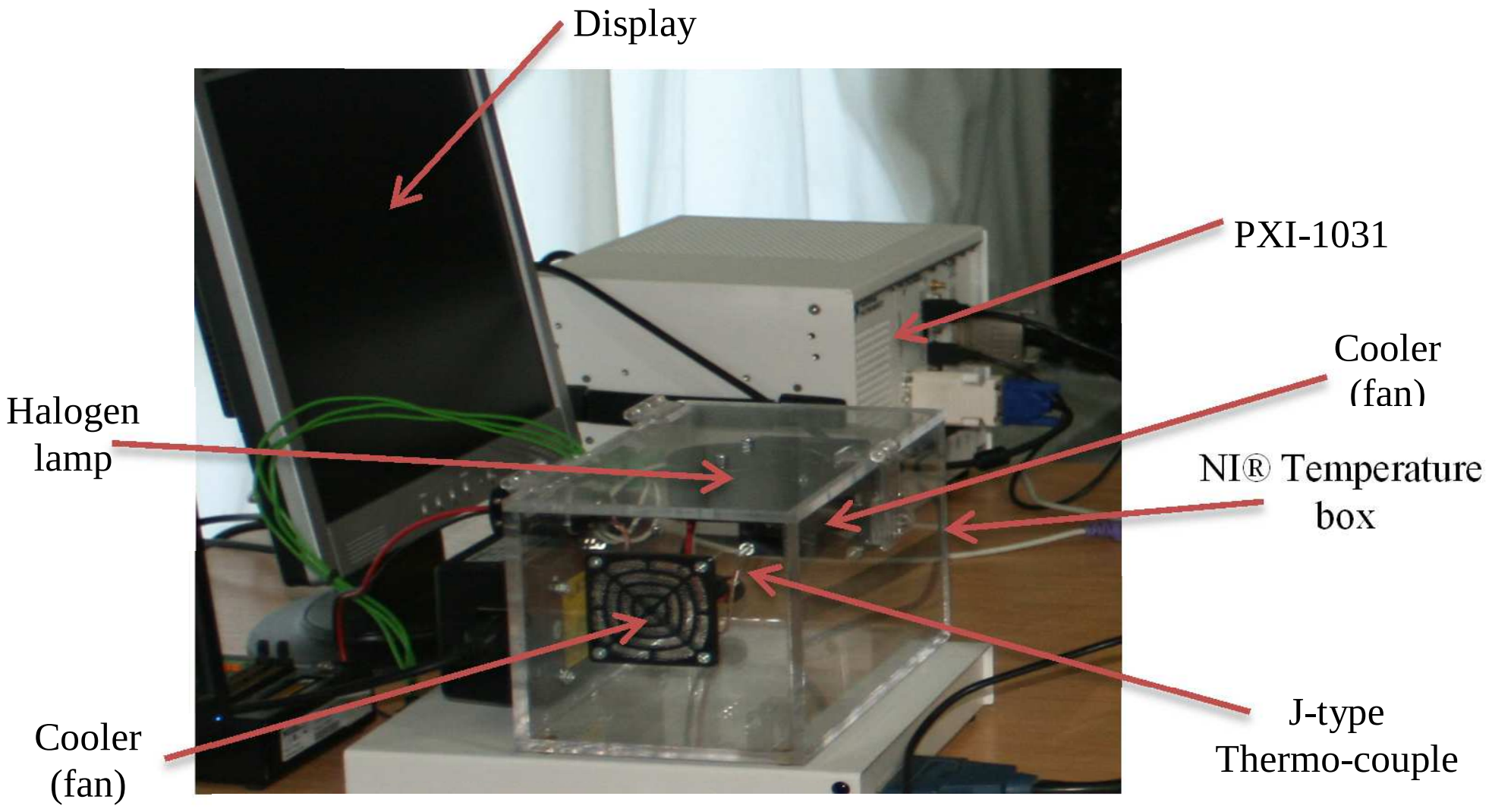

**Figure 1.** Experimental test bench.

We set the reference temperature to be higher than the ambient temperature using our input control model in NI® LabView environment and start collecting data from the J-type thermo-couple until the temperature reaches to our set point. The process duration depends on our set temperature value that is how much higher than the ambient temperature and where the thermo-couple is set, in other words, how far the thermo-couple from the heat source, viz. halogen lamp. We have chosen several different measurement points and temperature values, and data sampling values to collect test data inside the box as displayed in Fig. 2. In order to observe only temperature raise inside the box, a cooler (fan) of the system is switched-off and airflow window on the opposite side wall of the cooler is blocked. By this way, we have isolated the box from any additional energy flow-in and flow-out except for radiation heat in-flow that we have neglected in our studies because inside the laboratory room, only light source is luminescent lamps, which radiate extremely low heat energy.

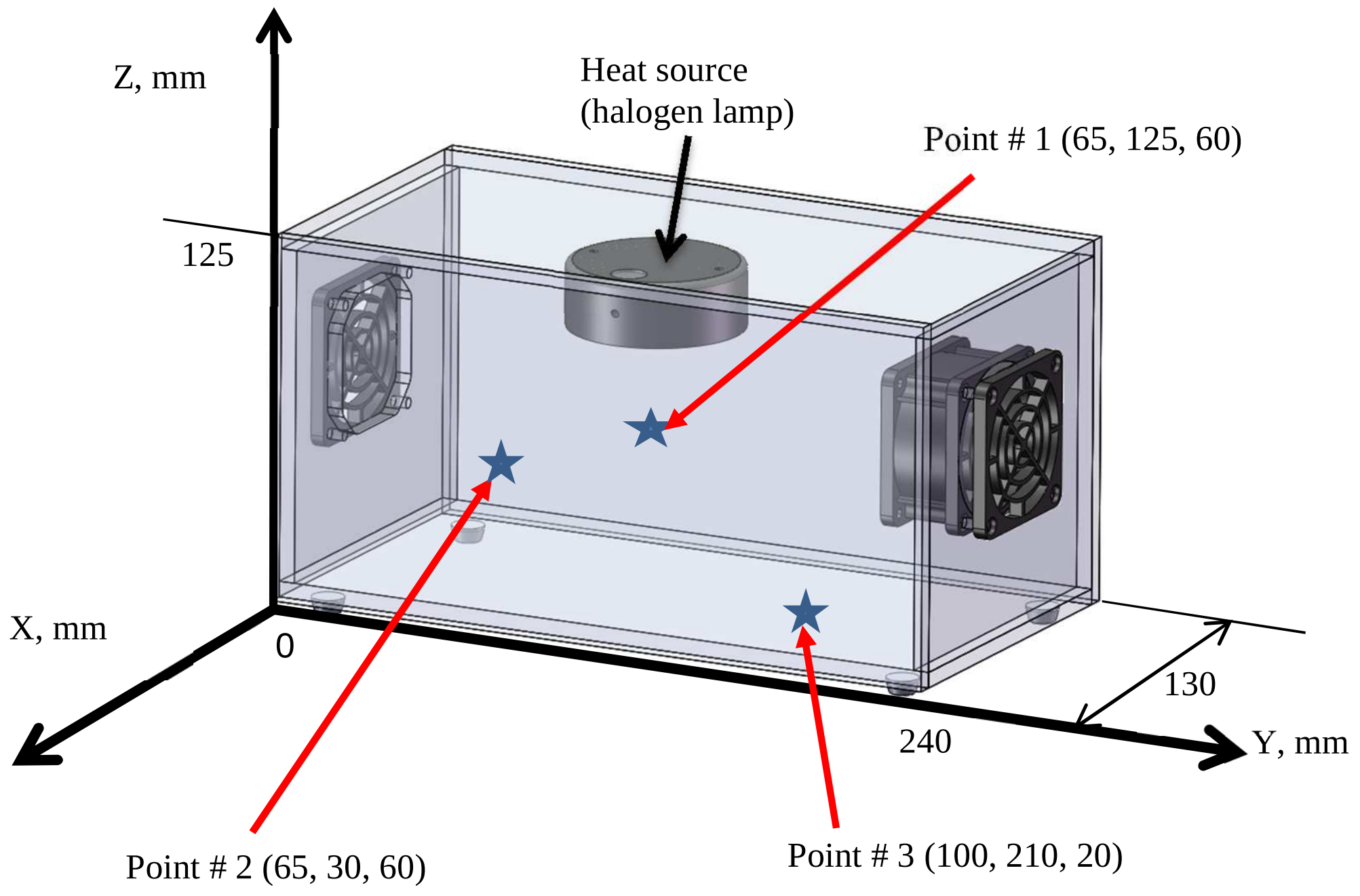

**Figure 2.** Experimental test data points.

## 4. Experimental data and their analysis

We have collected measurement data of temperature changes or heat flow rates from the J-type of thermo-couple inside the closed box from several points as shown in Fig. 2 with different data sampling frequencies (samples per second), such as 100 Hz, 500 Hz and 1000 Hz. For instance, from point 1, a number of data sets are collected with different sampling rates and different set temperatures as shown in Figure 3 and 4. From the trends of a heat flow process it is clear that the process of temperature raise evolves exponentially and there is a considerable noise present in the measured data. The level of a noise in the measured data is considerably larger particularly in higher sampling frequency (for 500 Hz and 1000 Hz) of our measured data sets. Therefore, we need to apply a digital filter to smooth the data before performing curve fitting.

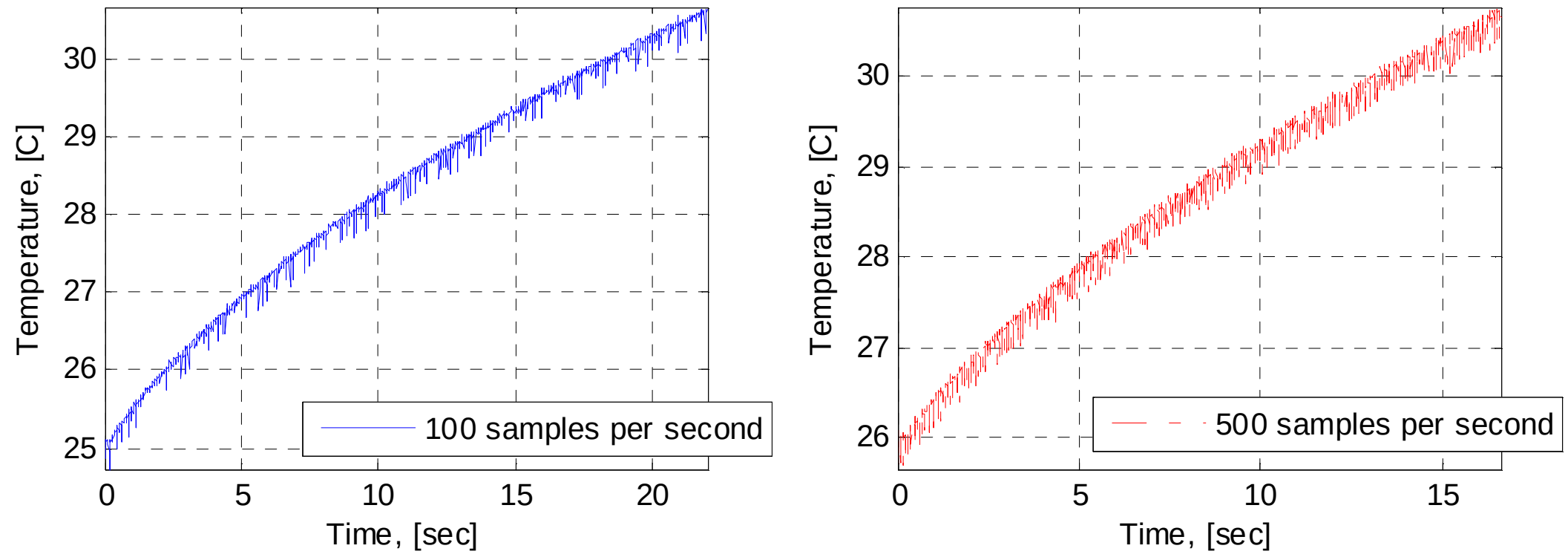

**Figure 3.** Experimental data collected in Point # 1 with 100 Hz and 500 Hz.

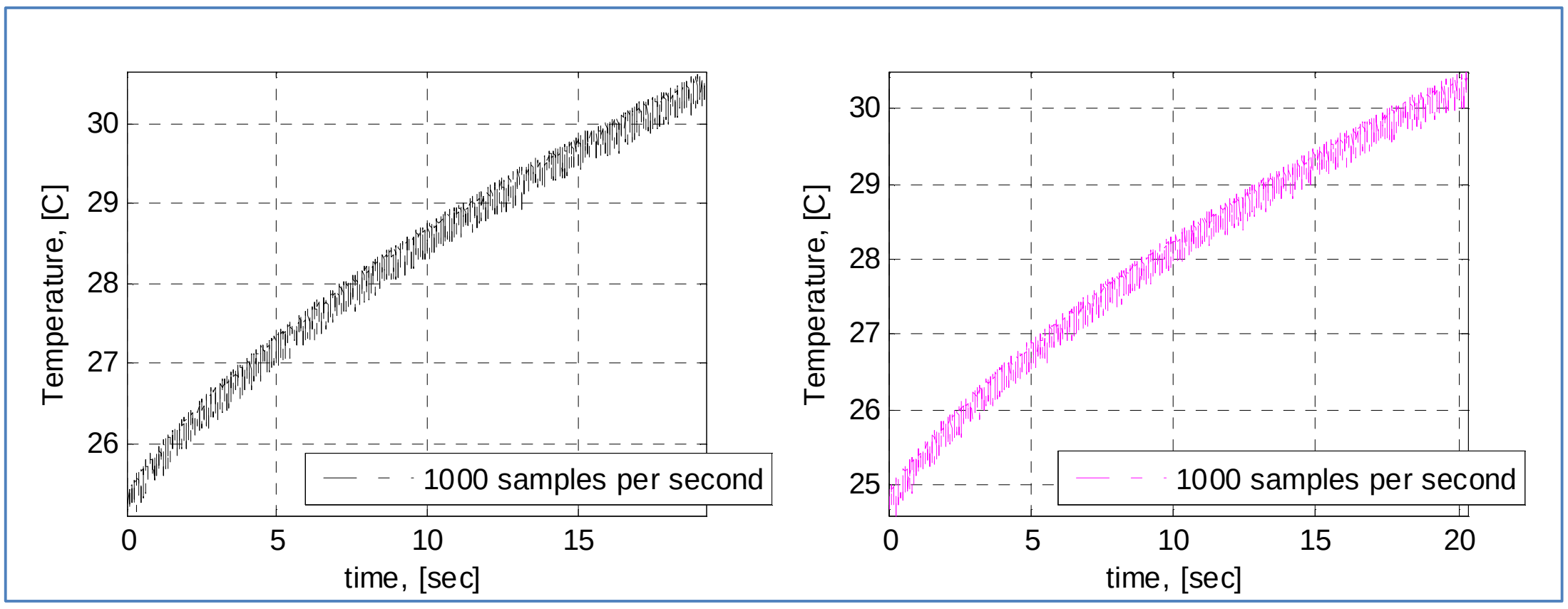

**Figure 4.** Experimental data collected in Point # 1 with 1000 Hz.

It is stated and demonstrated in several literature sources that Savitzky-Golay filter [9, 10] is a low-pass finite impulse response filter that is well suited to smooth noisy data. Besides, this type of filters are particularly well suited for time domain data sets. Therefore, we have employed the Savitzky-Golay filter to smooth our measured data from the J-type thermo-couple inside the temperature box. The filter has two parameters, viz. window size and order by which we can control quality of data smoothing. We have employed a MATLAB built-in function - sgolayfilt(DATA, Order, WinSize). We have used a third order filter with a window size of 901 that is smaller than half of our measurement data points.

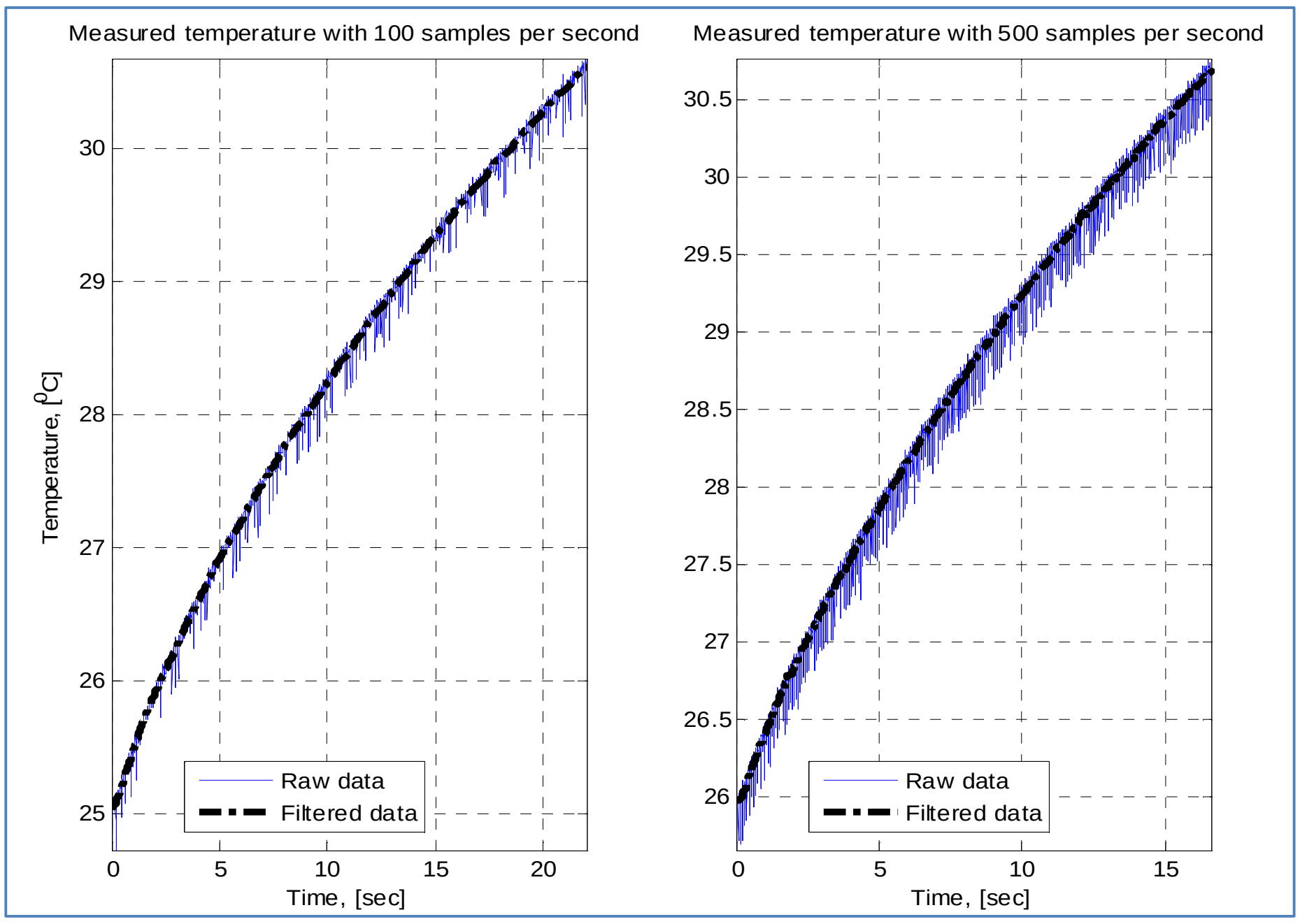

**Figure 5.** Raw experimental data (collected in Point # 1 with 100 Hz and 500 Hz) and its filtered part with the Savitzky-Golay filter.

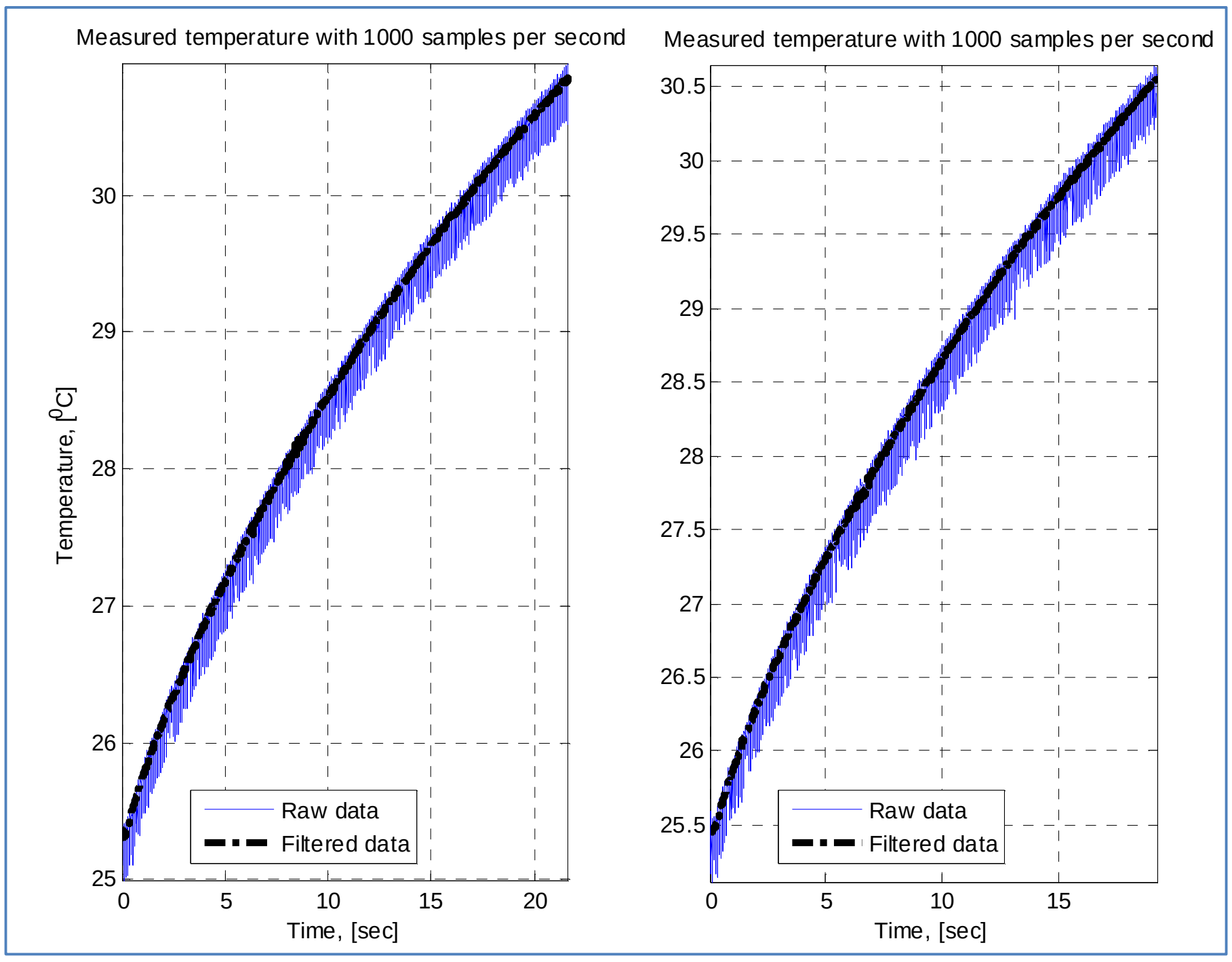

**Figure 6.** Experimental data (collected in Point # 1 with 1000 Hz) and its filtered component with the Savitzky-Golay filter.

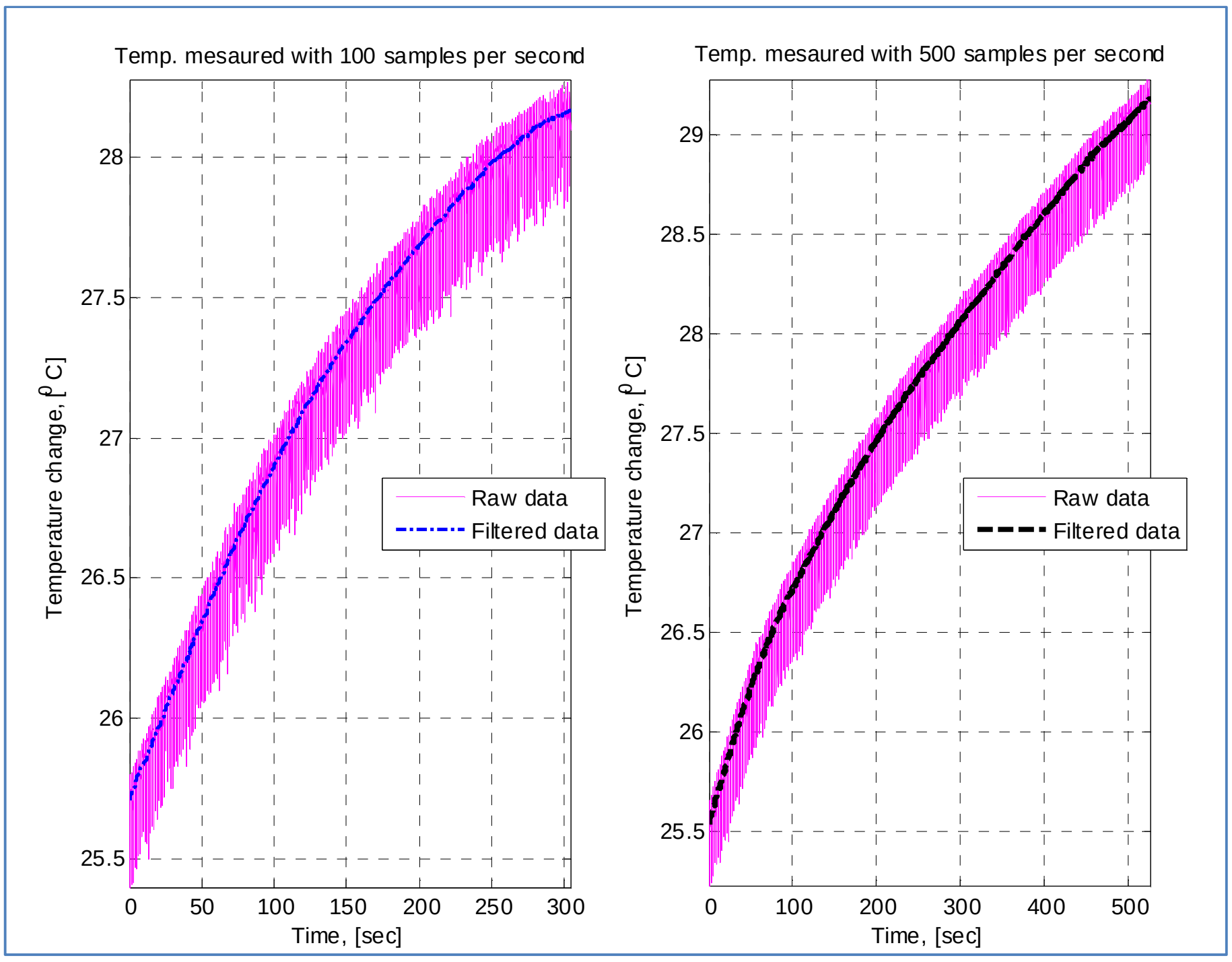

**Figure 7.** Raw experimental data (collected in Point # 2 with 100 Hz and 500 Hz) and its filtered part with the Savitzky-Golay filter.

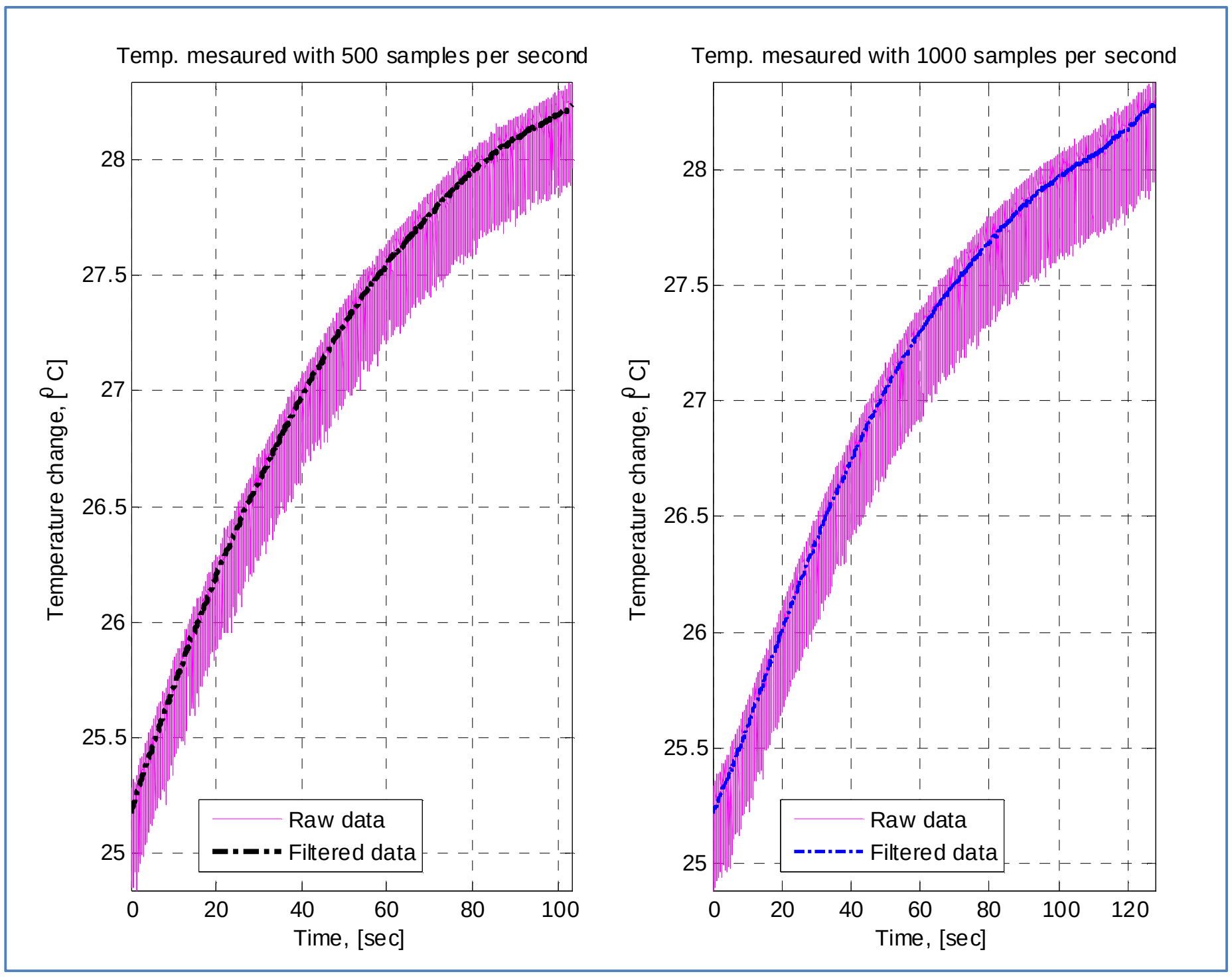

**Figure 8.** Raw experimental data (collected in Point # 3 with 500 Hz and 1000 Hz) and its filtered part with the Savitzky-Golay filter.

By observing the measured data and its filtered component with Savitzky-Golay filter displayed in Figure 3 – 6, it is clear that the filter has performed well and smoothed the raw data nicely and neatly. It is reasonably visible that the temperature change inside a closed box is in an exponential form. To compute a curve fit model of the measured data (smoothed ones) we first have applied the least squares method. The measured data are composed of two columns of data, one of which is time series and the other is a measured temperature change in ℃. By re-writing the expression in (13) and denoting the constant coefficients by: $a = \frac{T_a}{\tau}, b = K; c = \frac{1}{\tau}$, we get the system response formulation to be

$$f(t) = a * e^{-c*t} - b * e^{-c*t} + b \tag{14}$$

Where $t$ is the variable corresponding to time series. The formulation (14) is equivalent to (13) and in it there are three unknown parameters, viz. $a$, $b, c$ to be computed from the $n$ number of equations of the least squares method as formulated in (15) - (16). Our objective is to define a such formulation of $f_i(t)$ for which

$$\sum_{i=1}^{m} [f_i(t_i) - y_i]^2 \tag{15}$$

is minimized. Where $y_i$ is measured data points. In order to solve this minimization problem we need to compute such values of the unknowns that will minimize the expression (15).

$$F(a, b, c) = \sum_{i=1}^{m} [a * e^{-c*t_i} - b * e^{-c*t_i} + b - y_i]^2 \tag{16}$$

In order to minimize $F(a, b, c)$ with respect to $a, b, c$, we compute the partial derivatives and equate to zero:

$$\frac{\partial F(a, b, c)}{\partial a} = \frac{\partial F(a, b, c)}{\partial b} = \frac{\partial F(a, b, c)}{\partial c} = 0 \tag{17}$$

From (16), we shall have the following equations to be solved simultaneously:

$$\frac{\partial F(a,b,c)}{\partial a} = 2\sum_{i=1}^{m} e^{-c*t_i} \sum_{i=1}^{m} [a * e^{-c*t_i} - b * e^{-c*t_i} + b - y_i] \tag{18}$$

$$\frac{\partial F(a,b,c)}{\partial b} = \left(2\sum_{i=1}^{m} e^{-2*c*t_i} - 2\right) \sum_{i=1}^{m} [a * e^{-c*t_i} - b * e^{-c*t_i} + b - y_i] \tag{19}$$

$$\frac{\partial F(a,b,c)}{\partial c} = \left(-2a\sum_{i=1}^{m} t_i\, e^{-c*t_i} - 2b\sum_{i=1}^{m} t_i\, e^{-c*t_i}\right) \sum_{i=1}^{m} [a * e^{-c*t_i} - b * e^{-c*t_i} + b - y_i] \tag{20}$$

Solving the above shown three equations is not straight-forward due to nonlinearity of the equations. Thus, we shall employ the Levenberg-Marquardt method to solve these equations numerically. The Levenberg-Marquardt algorithm [11, 12] adaptively varies the parameter updates and the Marquardt's update relationship [12] is:

$$[J^T\, W\, J + \lambda\, diag(J^T W J)]h_{lm} = J^T\, W\, (y_i - \hat{y}(t_i; p)) \tag{21}$$

In the expression (21), $\hat{y}(t_i; p)$ is a fitting function of the parameters $p$ to a set of $m$ data points of $(t_i, y_i)$. The parameter update $h_{lm}$ with the Levenberg-Marquardt method moves the parameters towards the steepest descent. The Jacobian matrix $J$ (of the size of $m\ x\ n$) represents $\frac{\partial \hat{y}}{\partial p}$ that is the local sensitivity of the function $\hat{y}(t_i; p)$. The weighting matrix $W$ is diagonal with $1/\sigma_{yi}^2$ where $\sigma_{yi}$ is the measurement error of $y(t_i)$. The algorithmic parameter $\lambda$ defines how to take the update sizes; in other words, it influences on the iteration process in order to improve the solution search either by increasing it for worse approximation or decreasing it for better approximation (solutions).

There are many program codes created in Fortran, C and other programming languages, for instance, [13, 14, 15], and MATLAB's Curve Fitting toolbox™ to implement the Levenberg-Marquardt method for numerical solution search of a non-linear fit model's parameters. In our studies, we have compared the fit models computed from the two approaches MATLAB's Curve Fitting toolbox™ and M-files written in MATLAB based on the Levenberg-Marquardt method.

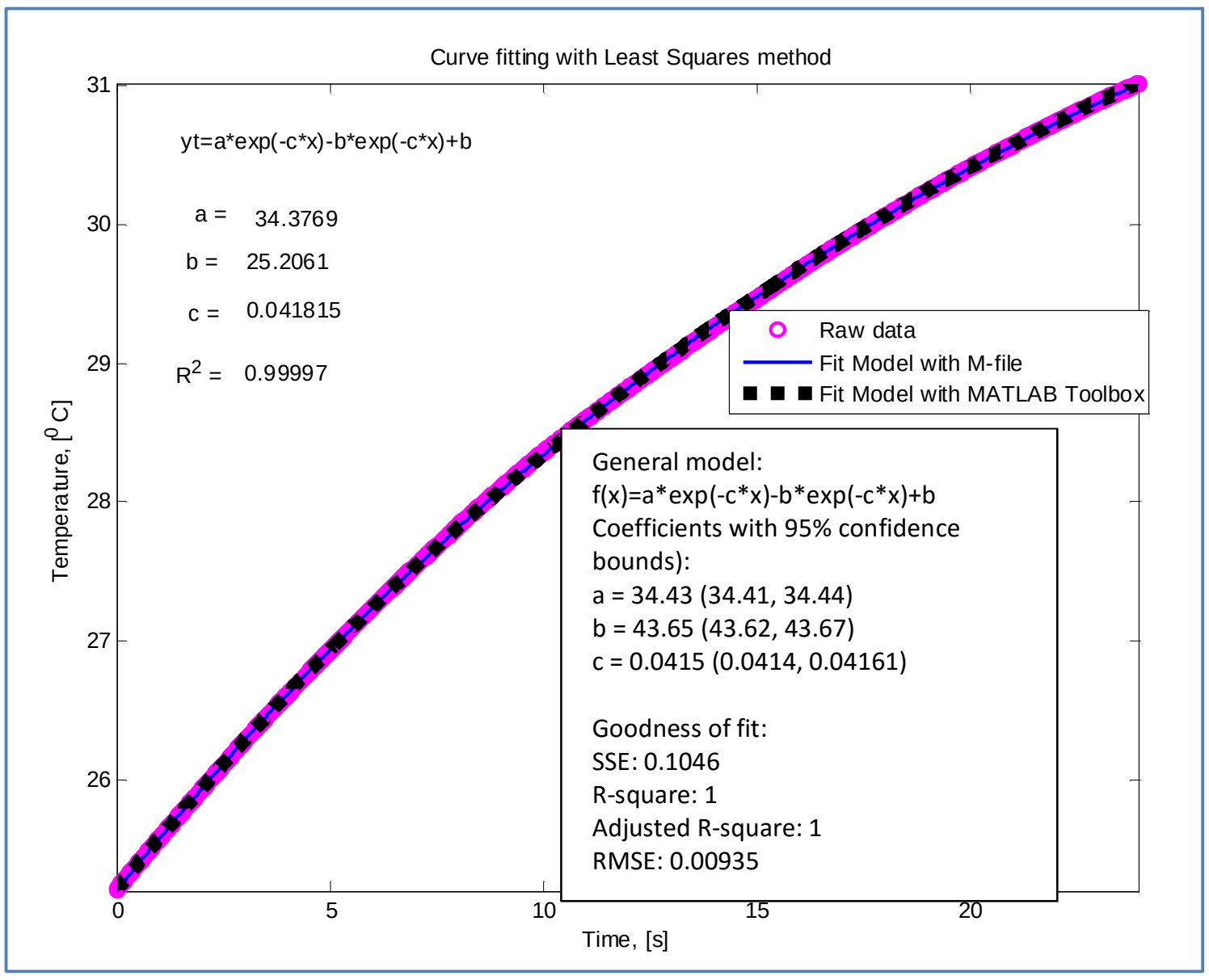

**Figure 9.** Raw experimental data (collected in Point # 1 with 100 Hz) vs. Fit model (M-file – upper left *yt* and Curve Fitting Toolbox – lower right *f(x)*).

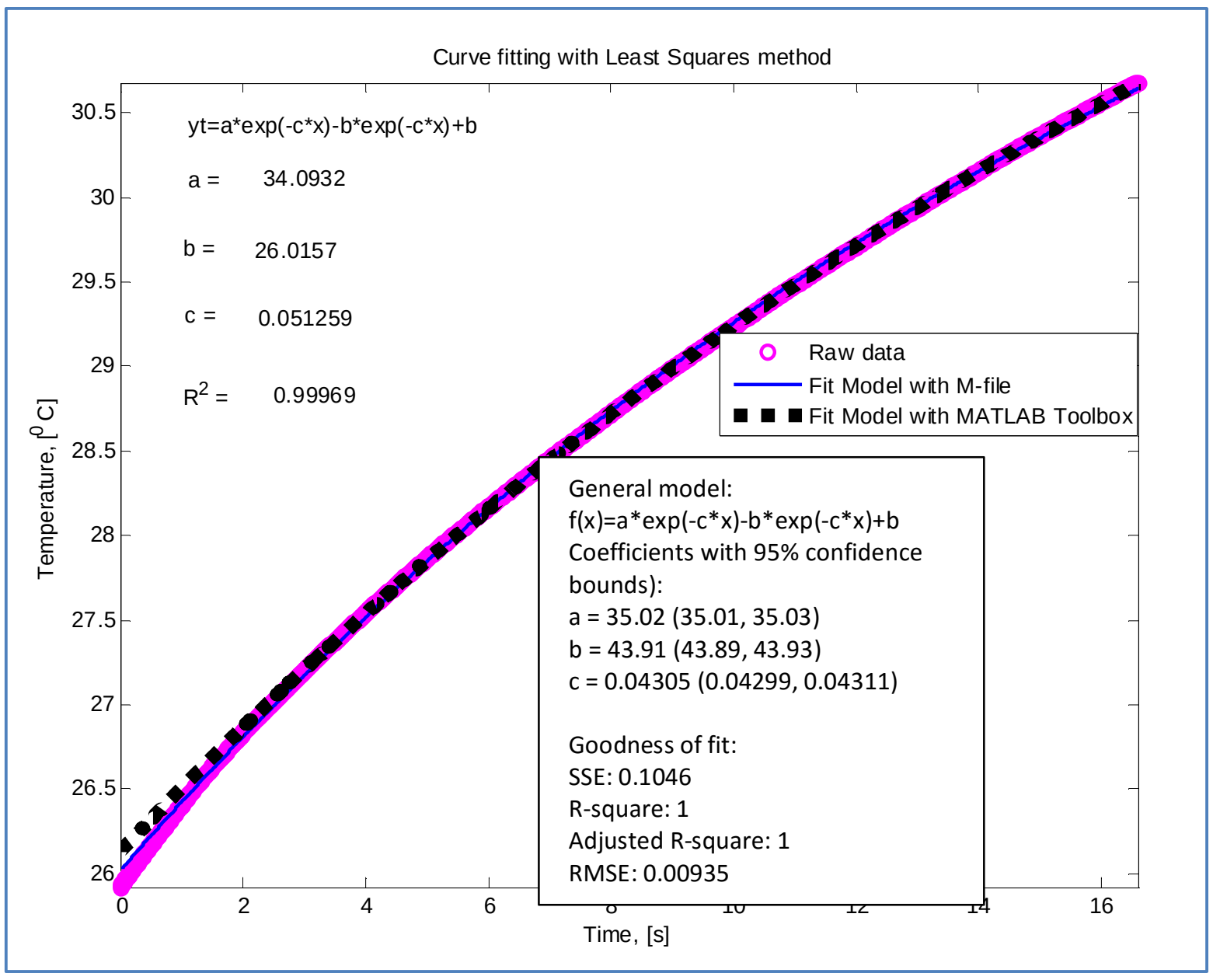

**Figure 10.** Raw experimental data (collected in Point # 1 with 500 Hz) vs. Fit model (M-file – upper left *yt* and Curve Fitting Toolbox – lower right *f(x)*).

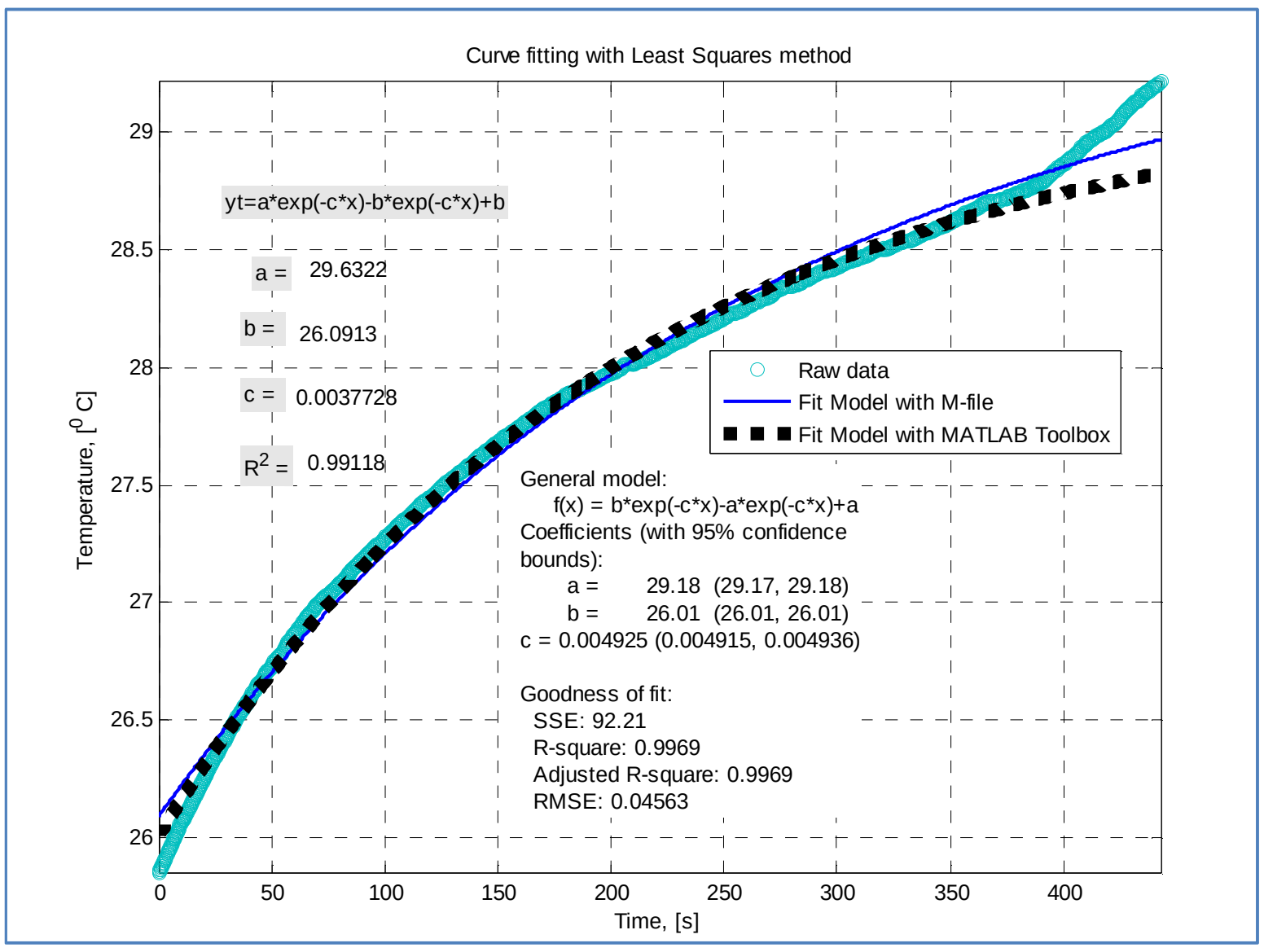

**Figure 11.** Raw experimental data (collected in Point # 2 with 100 Hz) vs. Fit model (M-file – upper left *yt* and Curve Fitting Toolbox – lower right *f(x)*).

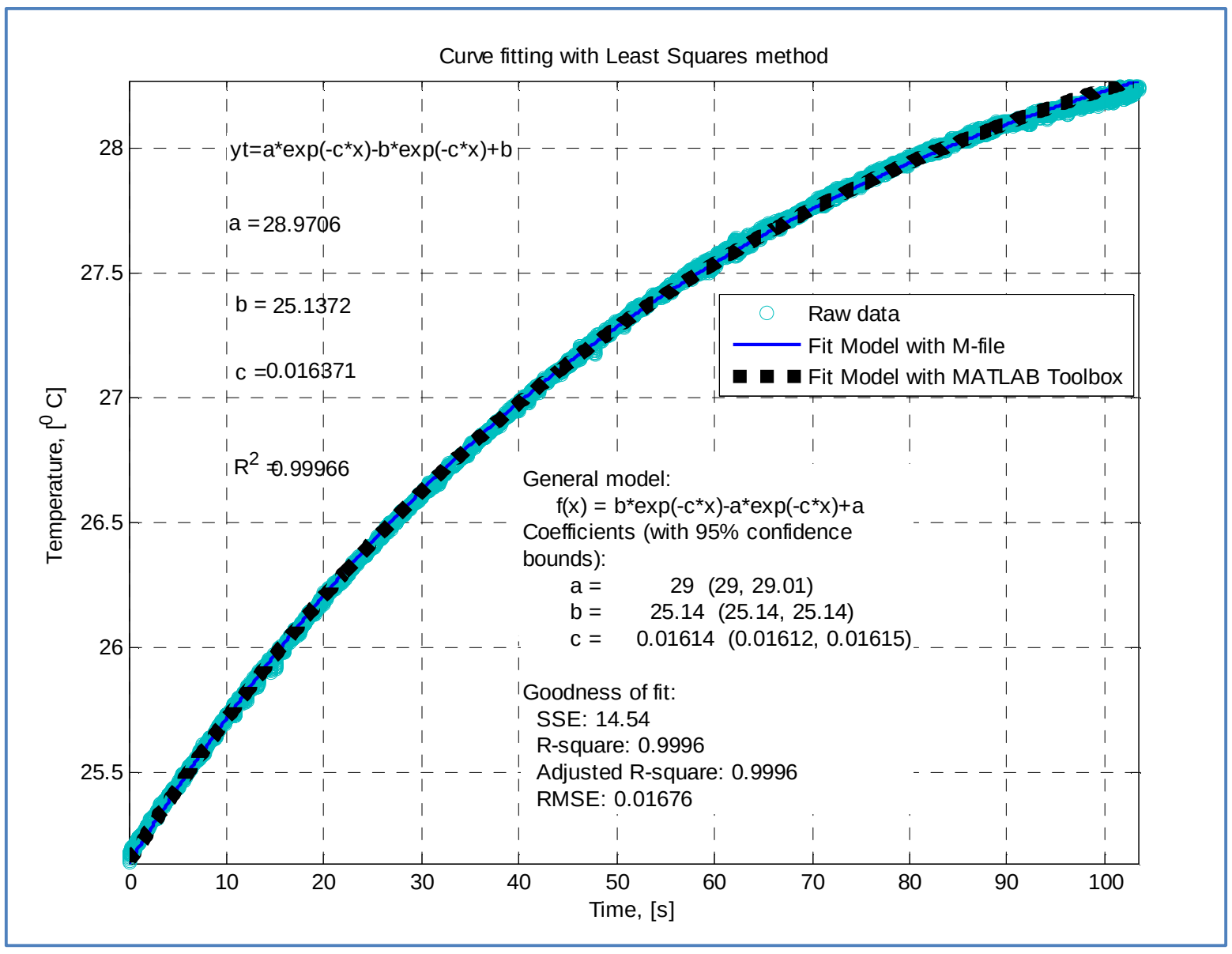

**Figure 12.** Raw experimental data (collected in Point # 2 with 500 Hz) vs. Fit model (M-file – upper left *yt* and Curve Fitting Toolbox – lower right *f(x)*).

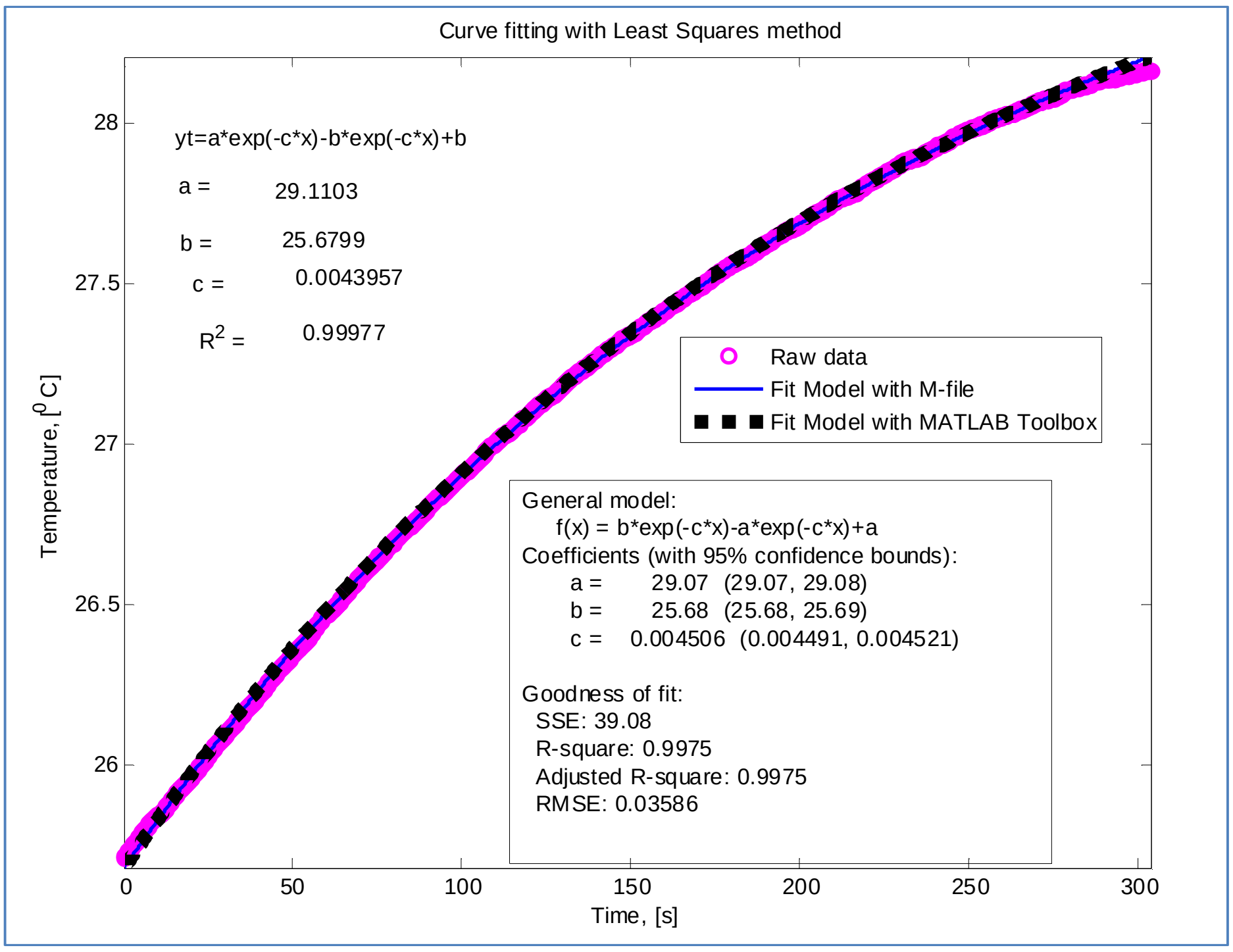

**Figure 13.** Raw data (100 samples in point #2) vs. fit models (M-file – upper left *yt* and Curve Fitting Toolbox – lower right *f(x)*).

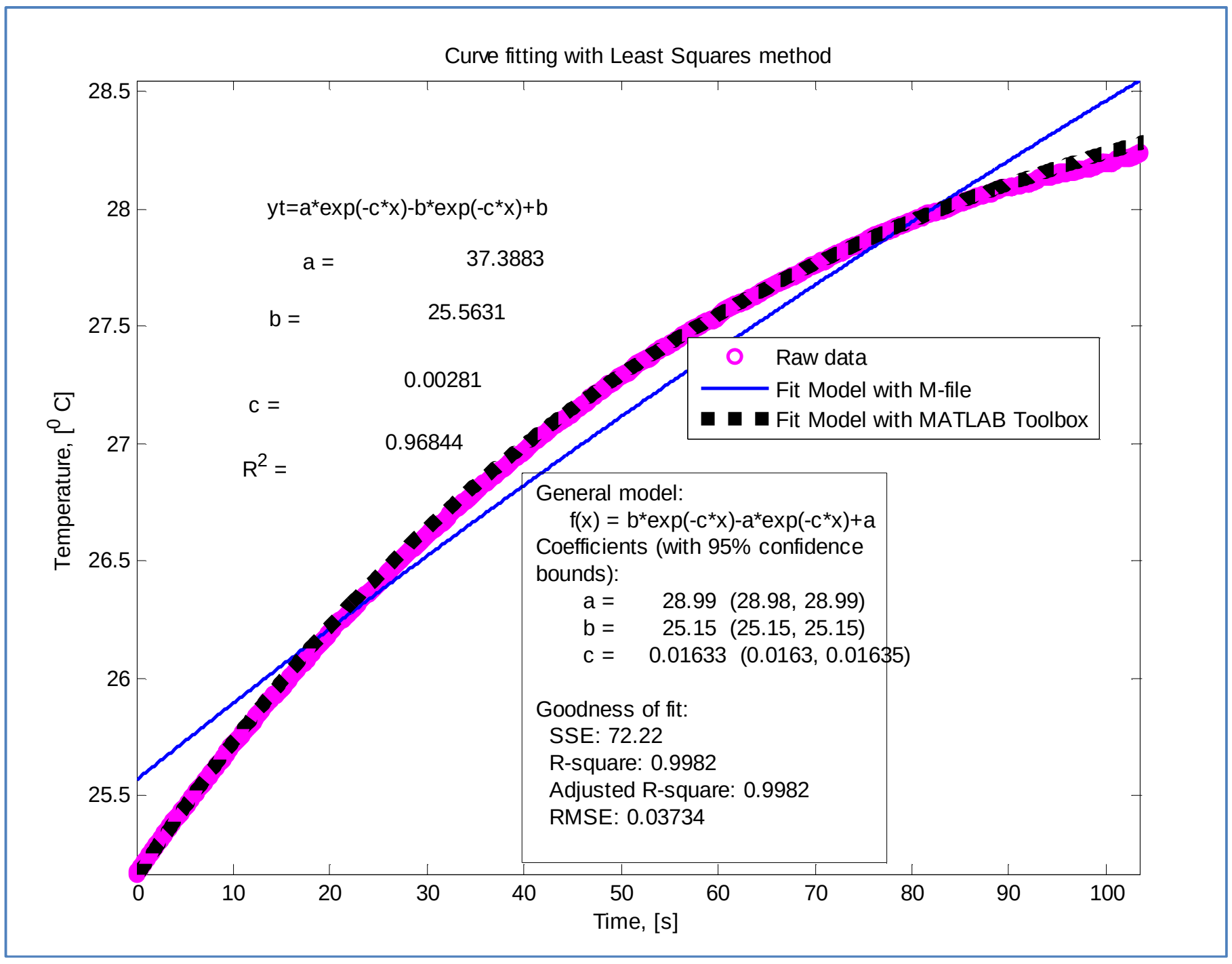

**Figure 14.** Raw data (500 samples per second at point #3) vs. fit models (M-file – upper left *yt* and Curve Fitting Toolbox – lower right *f(x)*).

The computed values of the parameters from the fit models are given in Table I below.

Table I. Parameter values defined from the curve fit models.

| Parameters | Point # 1 | | Point # 2 | | Point # 3 | |
|---|---|---|---|---|---|---|
| | MATLAB Curve fitting Toolbox | Code (M-file) | MATLAB Curve fitting Toolbox | Code (M-file) | MATLAB Curve fitting Toolbox | Code (M-file) |
| $a = \frac{T_a}{\tau}$ | 34.43 (100 Hz)<br>35.02 (500 Hz) | 34.37 (100 Hz)<br>34.09 (500 Hz) | 29.18(100Hz)<br>29 (500 Hz) | 29.63(100Hz)<br>26.09(500Hz) | 29.07(100Hz)<br>28.99(500Hz) | 29.11(100Hz)<br>37.39(500Hz) |
| b = K | 43.65 (100 Hz)<br>43.91 (500 Hz) | 25.21 (100 Hz)<br>26.01 (500 Hz) | 26.01(100Hz)<br>25.14(500Hz) | 28.97(100Hz)<br>25.14(500Hz) | 25.68(100Hz)<br>25.15(500Hz) | 25.68(100Hz)<br>25.56(500Hz) |
| $c = \frac{1}{\tau}$ | 0.0415(100Hz)<br>0.043 (500 Hz) | 0.0418(100Hz)<br>0.0512(500Hz) | 0.0049(100Hz)<br>0.0161(500Hz) | 0.004(100Hz)<br>0.016(500Hz) | 0.004(100Hz)<br>0.016(500Hz) | 0.004(100Hz)<br>0.002(500Hz) |

All computed curve fit models of the measured data from Point # 1, 2 and 3 with 100Hz and 500 Hz sampling frequencies with fit models are shown in Fig. 9-14 by using the least squares method from MATLAB Curve fitting toolbox and M-file scripts with the Levenberg-Marquardt algorithm are well converged and coefficient of determination ($R^2$) is very close to 1.0. The two fit model coefficients found from the Curve Fitting Toolbox and M-files are close and the parameter values for $a = \frac{T_a}{\tau}$, b = K, and $c = \frac{1}{\tau}$ are to be considerable stable for point # 2 and 3. Fit models found with the M-files with the Levenberg-Marquardt algorithm are considerably nice and well converged with the Toolbox. Different sampling data values have not made substantial influence on the values of the fit models. Besides, different starting initial temperature (ambient temperature) values have not influenced either on the found parameters. On the other hand, some small differences in parameter (a, b, c) values can be explained with the linear property of J-type thermo-couples' thermoelectric voltage as a function of temperature (°C) that is shown in Fig. 15.

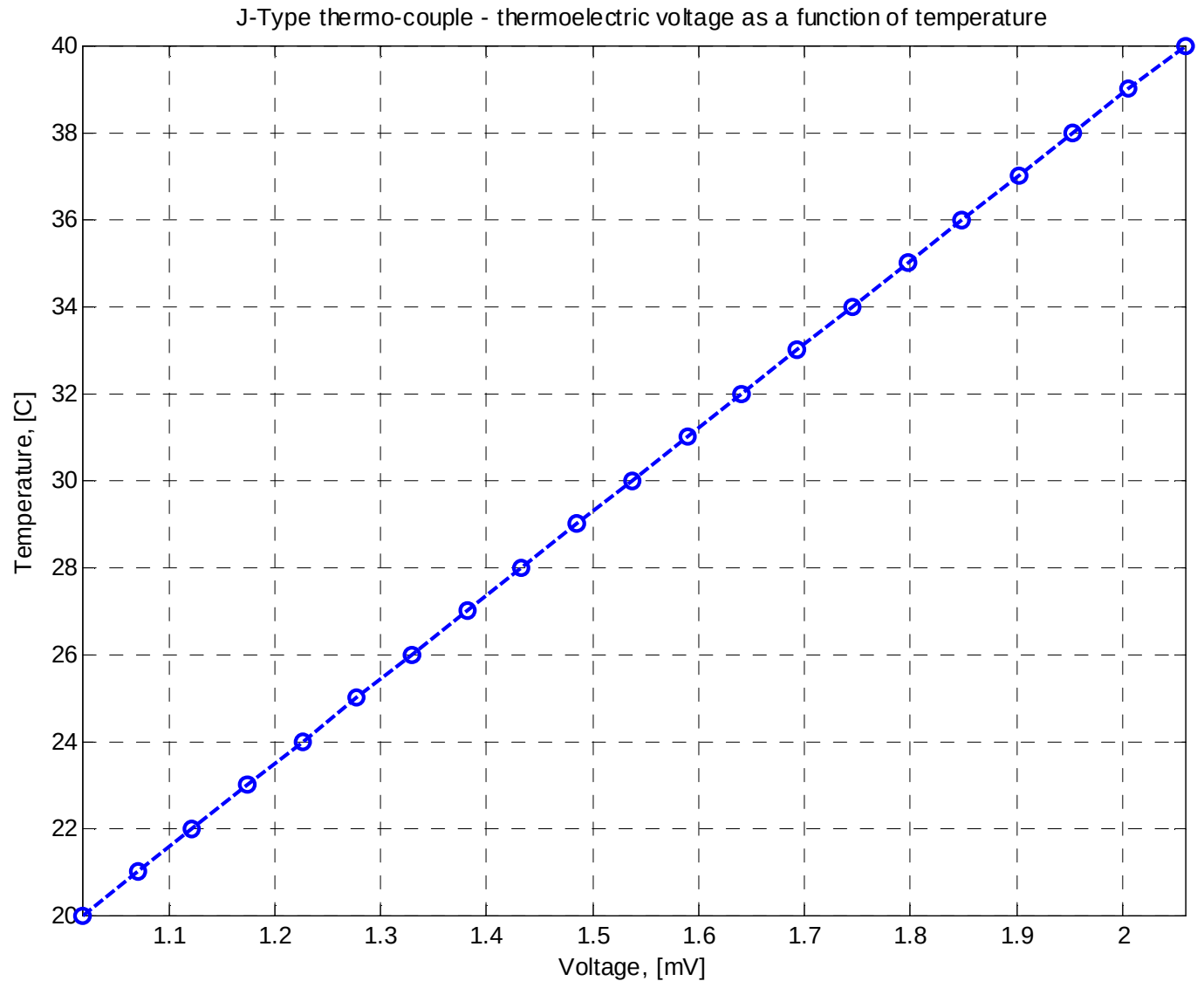

**Figure 15.** J-type thermo-couple reference junctions at 20 … 40°C.

## Conclusions and Discussions

We have studied temperature/heat flow process inside a closed environment with the test rig of the NI Temperature box with respect to our derived mathematical formulations of the process. The test data measured from the test object have showed that the derived process model describe the process itself adequately well and the computed parameter values of the process model with the least squares method by using the Levenberg-Marquardt method are considerably stable for different data sampling and measurement points inside the box. Our further studies will be dedicated to designing control algorithm and controller of the heat flow process around the reference temperature value inside a box that has heat source and heat outflow (blower/fan).

**Acknowledgements.** Part of this research study is developed within the EC TEMPUS grant 543-TEMPUS-1-2013-1-SE-JPCR.